\documentclass[final,leqno]{siamltex}

\usepackage[latin1]{inputenc}
\usepackage{psfrag}

\usepackage{graphicx}
\usepackage{amssymb,amsfonts,cmmib57}

\newcommand{\myfigure}[2]{\includegraphics[#1]{#2}}

\newcommand{\mcgq}{$\mathrm{mcG}(q)$}
\newcommand{\mdgq}{$\mathrm{mdG}(q)$}

\newcommand{\PD}[2]{\frac{\partial#1}{\partial#2}}
\newcommand{\real}{\mathbb{R}}

\newtheorem{remark}{{\it Remark}}[section]

\title{Multi-Adaptive Galerkin Methods for ODEs
II:\\ implementation and Applications\thanks{Received by the
editors May 23, 2001; accepted for publication (in revised form)
May 1, 2003; published electronically December 5, 2003.
\URL sisc/25-4/38973.html}}

\author{Anders Logg\thanks{Department of Computational Mathematics,
Chalmers University of Technology, SE--412 96 G\"{o}teborg, Sweden
(logg@math.chalmers.se).}}

\begin{document}
\maketitle
\vspace{-1.2in}
\slugger{sisc}{2003}{25}{4}{1119--1141}
\vspace{.9in}

\setcounter{page}{1119}

\begin{abstract}
Continuing the discussion of the multi-adaptive Galerkin methods
{$\mathrm{mcG}(q)$} and {$\mathrm{mdG}(q)$} presented in [A. Logg, {\it SIAM J. Sci. Comput.}, 24 (2003),
pp.~1879--1902],
we present adaptive algorithms for global error control,
iterative solution methods for the discrete equations,
features of the implementation \emph{Tanganyika}, and
computational results for a variety of ODEs. Examples include the
Lorenz system, the solar system, and a number of time-dependent
PDEs.
\end{abstract}



\pagestyle{myheadings}
\thispagestyle{plain}
\markboth{ANDERS LOGG}{MULTI-ADAPTIVE GALERKIN METHODS FOR ODEs
II}
\setlength\parskip{1\parskip}
\section{Introduction}

In this paper we apply the multi-adaptive Galerkin methods
\mcgq\ and \mdgq, presented in
\cite{logg:multiadaptivity:I}, to a variety of problems chosen to illustrate
the potential of multi-adaptivity. Throughout this paper, we
solve the ODE initial value problem
\begin{equation}
\left\{
\begin{array}{rcl}
\dot{u}(t) &=& f(u(t),t),\ t\in(0,T], \\
u(0) &=& u_0,
\end{array}
\right.
\label{eq:u'=f}
\end{equation}
where $u : [0,T] \rightarrow \real^N$, $f : \real^N \times (0,T]
\rightarrow \real^N$ is a given bounded function that is
Lipschitz continuous in $u$, $u_0 \in \real^N$ is a given initial
condition, and $T>0$ a given final time.

We refer to \cite{logg:multiadaptivity:I} for a detailed
description of the multi-adaptive methods. Here we recall that
each component $U_i(t)$ of the approximate solution $U(t)$ is a
piecewise polynomial of degree $q_i=q_i(t)$ on a partition of
$(0,T]$ into $M_i$ subintervals of lengths
$k_{ij}=t_{ij}-t_{i,j-1}$, $j=1,\ldots,M_i$. On the interval
$I_{ij}=(t_{i,j-1},t_{ij}]$, component $U_i(t)$ is thus a
polynomial of degree $q_{ij}$.

Before presenting the examples, we discuss adaptive algorithms for
global error control and iterative solution methods for the
discrete equations. We also give a short description of the
implementation \emph{Tanganyika}.

\section{Adaptivity}

In this section we describe how to use the a posteriori
error estimates presented in \cite{logg:multiadaptivity:I} in an
adaptive algorithm.

\subsection{A strategy for adaptive error control}

The goal of the algorithm is to produce an approximate solution
$U(t)$ to (\ref{eq:u'=f}) within a given tolerance $\mathrm{TOL}$
for the error $e(t)=U(t)-u(t)$ in a given norm $\| \cdot \|$. The
adaptive algorithm is based on the a posteriori error estimates,
presented in \cite{logg:multiadaptivity:I}, of the form
\begin{equation}
\| e \| \leq \sum_{i=1}^N\sum_{j=1}^{M_i}
k_{ij}^{p_{ij}+1} r_{ij} s_{ij}
\label{eq:estimate,1}
\end{equation}
or
\begin{equation}
\| e \| \leq \sum_{i=1}^N
S_i \max_{j} k_{ij}^{p_{ij}} r_{ij},
\label{eq:estimate,2}
\end{equation}
where $\{s_{ij}\}$ are $\emph{stability weights}$,
$\{S_i\}$ are $\emph{stability factors}$ (including
interpolation constants), $r_{ij}$ is a local measure of the
residual $R_i(U,\cdot)=\dot{U}_i - f(U,\cdot)$ of the approximate
solution $U(t)$, and where we have $p_{ij}=q_{ij}$ for \mcgq\ and
$p_{ij}=q_{ij}+1$ for \mdgq.

We use (\ref{eq:estimate,2}) to determine the individual
time-steps, which should then be chosen as
\begin{equation}
k_{ij} = \left( \frac{\mathrm{TOL}/N}{S_i \ r_{ij}}
\right)^{1/{p_{ij}}}.
\label{eq:timestep}
\end{equation}
We use (\ref{eq:estimate,1}) to evaluate the resulting error at
the end of the computation, noting that (\ref{eq:estimate,1}) is
sharper than (\ref{eq:estimate,2}).

The adaptive algorithm may then be expressed as follows: Given a
tolerance $\mathrm{TOL}>0$, make a preliminary guess for the
stability factors and then
\begin{enumerate}
\item
solve the primal problem with time-steps based on
(\ref{eq:timestep}).
\item
solve the dual problem and compute stability factors and
stability weights.
\item
compute an error estimate $E$ based on (\ref{eq:estimate,1}).
\item
if $E \leq \mathrm{TOL}$, then stop, and if not, go back to (i).
\end{enumerate}

Although this seems simple enough, there are some difficulties
involved. For one thing, choosing the time-steps based on
(\ref{eq:timestep}) may be difficult, since the residual depends
implicitly on the time-step. Furthermore, we have to choose the
proper data for the dual problem to obtain a meaningful error
estimate. We now discuss these issues.

\subsection{Regulating the time-step}

To avoid the implicit dependence on $k_{ij}$ for $r_{ij}$ in
(\ref{eq:timestep}), we may try replacing (\ref{eq:timestep}) with
\begin{equation}
k_{ij} = \left( \frac{\mathrm{TOL}/N}{S_i \ r_{i,j-1}} \right)^{1/{p_{ij}}}.
\label{eq:timestep,alternative}
\end{equation}
Following this strategy, if the time-step on an interval is small
(and thus also is the residual), the time-step for the next interval
will be large, so that (\ref{eq:timestep,alternative}) introduces
unwanted oscillations in the size of the time-step. We therefore
try to be a bit more conservative when choosing the time-step to
obtain a smoother time-step sequence. For
(\ref{eq:timestep,alternative}) to work, the time-steps on
adjacent intervals need to be approximately the same, and so we
may think of choosing the new time-step as the (geometric) mean
value of the previous time-step and the time-step given by
(\ref{eq:timestep,alternative}). This works surprisingly well for
many problems, meaning that the resulting time-step sequences are
comparable to what can be obtained with more advanced regulators.

We have also used standard \emph{PID} (or just \emph{PI})
regulators from control theory with the goal of satisfying
\begin{equation}
S_i k_{ij}^{p_{ij}}r_{ij} = \mathrm{TOL} / N,
\label{eq:regler,1}
\end{equation}
or, taking the logarithm with $C_i = \log(\mathrm{TOL}/(N S_i))$,
\begin{equation}
p_{ij} \log k_{ij} + \log r_{ij} = C_i,
\label{eq:regler,2}
\end{equation}
with maximal time-steps $\{k_{ij}\}$, following work by
S\"{o}derlind \cite{soderlind:control} and Gustafsson, Lundh, and
S\"{o}derlind \cite{soderlind:pid}.  This type of regulator
performs a little better than the simple approach described
above, provided the parameters of the regulator are well tuned.

\subsection{Choosing data for the dual}

Different choices of data $\varphi_T$ and $g$ for the dual
problem give different error estimates, as described in
\cite{logg:multiadaptivity:I}, where estimates for the quantity
\begin{displaymath}
L_{\varphi_T,g}(e) =
(e(T),\varphi_T) + \int_0^T (e,g) \ dt
\end{displaymath}
were derived. The simplest choices are $g=0$ and
$(\varphi_T)_i=\delta_{in}$ for control of the final time error
of the $n$th component. For control of the $l^2$-norm of the
error at final time, we take $g=0$ and $\varphi_T =
\tilde{e}(T)/\|\tilde{e}(T)\|$ with an approximation $\tilde{e}$
of the error $e$. Another possibility is to take $\varphi_T=0$
and $g_i(t) = \delta_{in}$ for control of the average error in
component $n$.

If the data for the dual problem are incorrect, the error
estimate may also be incorrect: with $\varphi_T$ (or $g$)
orthogonal to the error, the error representation gives only
$0\leq \mathrm{TOL}$. In practice, however, the dual---or at
least the stability factors---seems to be quite insensitive to
the choice of data for many problems so that it is, in fact,
possible to guess the data for the dual.

\subsection{Adaptive quadrature}

In practice, integrals included in the formulation of the two methods \mcgq\ and
\mdgq\ have to be evaluated using numerical quadrature.
To control the resulting quadrature error, the quadrature rule can
be chosen adaptively, based on estimates of the quadrature error
presented in \cite{logg:multiadaptivity:I}.

\subsection{Adaptive order, \boldmath$q$-adaptivity}

The formulations of the methods include the possibility of
individual and varying orders for the different components, as
well as different and varying time-steps. The method is thus
$q$\emph{-adaptive} (or $p$-adaptive) in the sense that the order
can be changed. At this stage, however, lacking a strategy for
when to increase the order and when to decrease the time-step,
the polynomial orders have to be chosen in an ad hoc
fashion for every interval. One way to choose time-steps and
orders could be to solve over some short time-interval with
different time-steps and orders, and optimize the choice of
time-steps and orders with respect to the computational time
required for achieving a certain accuracy. If we suspect that the
problem will change character, we will have to repeat this
procedure at a number of control points.


\section{Solving the discrete equations}
\label{sec:solving}

In this section we discuss how to solve the discrete equations that we
obtain when we discretize (\ref{eq:u'=f}) with the multi-adaptive
Galerkin methods. We do this in two steps. First, we present a
simple explicit strategy, and then we extend this strategy to an
iterative method.

\subsection{A simple strategy}

As discussed in \cite{logg:multiadaptivity:I}, the nonlinear
discrete algebraic equations for the \mcgq\ method (including
numerical quadrature) to be solved on every local interval
$I_{ij}$ take the form
\begin{equation}
\xi_{ijm} =
\xi_{ij0} +
k_{ij} \sum_{n=0}^{q_{ij}} w_{mn}^{[q_{ij}]} \ f_i(U(\tau_{ij}^{-1}(s_{n}^{[q_{ij}]})),\tau_{ij}^{-1}(s_n^{[q_{ij}]})), \ m = 1,\ldots,q_{ij},
\label{eq:equations}
\end{equation}
where $\{\xi_{ijm}\}_{m=1}^{q_{ij}}$ are the degrees of freedom
to be determined for component $U_i(t)$ on interval $I_{ij}$,
$\{w_{mn}^{[q_{ij}]}\}_{m=1,n=0}^{q_{ij}}$ are weights,
$\{s_n^{[q_{ij}]}\}_{n=0}^{q_{ij}}$ are quadrature points, and
$\tau_{ij}$ maps $I_{ij}$ to $(0,1]$: $\tau_{ij}(t) =
(t-t_{i,j-1})/(t_{ij}-t_{i,j-1})$. The discrete equations for the
\mdgq\ method are similar in structure and so we focus on the
\mcgq\ method.

The equations are conveniently written in fixed point form, so we
may apply fixed point iteration directly to (\ref{eq:equations});
i.e., we make an initial guess for the values of
$\{\xi_{ijm}\}_{m=1}^{q_{ij}}$, e.g., $\xi_{ijm} = \xi_{ij0}$ for
$m=1,\ldots,q_{ij}$, and then compute new values for these
coefficients from (\ref{eq:equations}), repeating the procedure
until convergence.

Note that component $U_i(t)$ is coupled to all other
components through the right-hand side $f_i=f_i(U,\cdot)$. This
means that we have to know the solution for all other components
in order to compute $U_i(t)$. Conversely, we have to know
$U_i(t)$ to compute the solutions for all other components, and
since all other components step with different time-steps, it
seems at first very difficult to solve the discrete equations
(\ref{eq:equations}).

As an initial simple strategy we may try to solve the
system of nonlinear equations (\ref{eq:equations}) by direct
fixed point iteration. All unknown values, for the component
itself and all other \emph{needed} components, are interpolated
or extrapolated from their latest known values. Thus, if for
component $i$ we need to know the value of component $l$ at some
time $t_i$, and we only know values for component $l$ up to
time $t_l < t_i$, the strategy is to extrapolate $U_l(t)$ from
the interval containing $t_l$ to time $t_i$, according to the
order of $U_l(t)$ on that interval.

\begin{figure}[h]
\begin{center}
\leavevmode
\psfrag{0}{$0$}
\psfrag{T}{$T$}
\includegraphics[width=\textwidth,height=8cm]{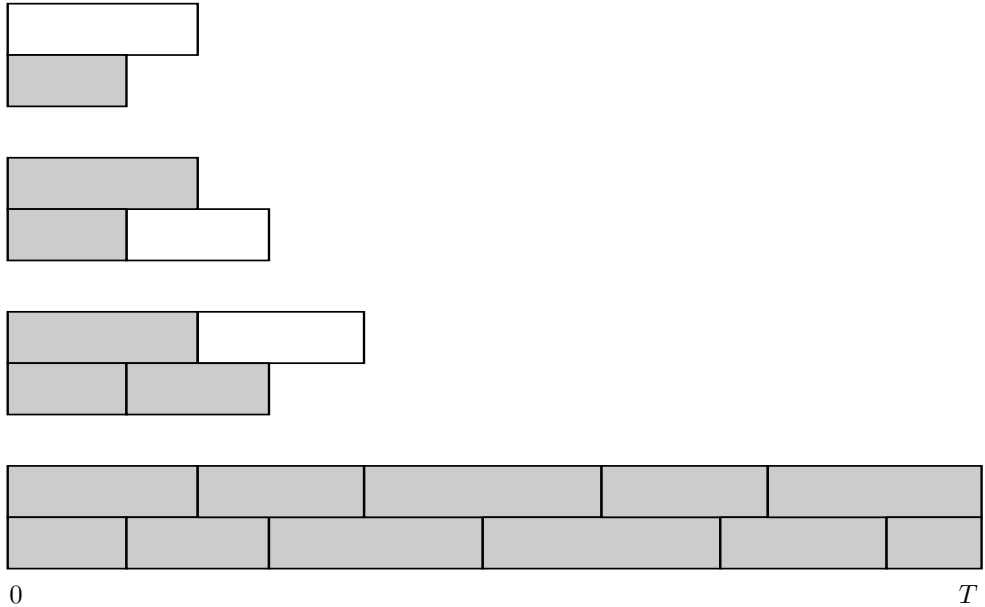}
\caption{The last component steps first and all needed
values are extrapolated or interpolated.}
\label{fig:stepping}
\end{center}
\end{figure}

In what order should the components now make their steps?
Clearly, to update a certain component on a specific interval, we
would like to use the best possible values of the other
components. This naturally leads to the following strategy:
\begin{equation}
\mbox{\textit{The last component steps first.}}
\label{eq:strategy}
\end{equation}
This means that we should always make a step with the component
that is closest to time $t=0$. Eventually (or after one
step), this component catches up with one of the other components,
which then in turn becomes the last component, and the procedure
continues according to the strategy (\ref{eq:strategy}), as
described in Figure \ref{fig:stepping}.

This gives an explicit time-stepping method in which
each component is updated individually once, following
(\ref{eq:strategy}), and in which we never go back to correct mistakes.
This corresponds to fixed point iterating once on the discrete
equations (\ref{eq:equations}), which implicitly define the
solution. We now describe how to extend this explicit
time-stepping strategy to an iterative process, in which the
discrete equations are solved to within a prescribed tolerance.

\subsection{An iterative method}

To extend the explicit strategy described in the previous section
to an iterative method, we need to be able to go back and redo
iterations if necessary. We do this by arranging
the elements---we think of an element as a component $U_i(t)$
on a local interval $I_{ij}$---in a \emph{time-slab}. This
contains a number of elements, a minimum of $N$ elements, and
moves forward in time. On every time-slab, we have to solve a
large system of equations, namely, the system composed of the
element equations (\ref{eq:equations}) for every element within
the time-slab. We solve this system of equations iteratively, by
direct fixed point iteration, or by some other method as
described below, starting from the last element in the time-slab,
i.e., the one closest to $t=0$, and continuing forward to the
first element in the time-slab, i.e., the one closest to $t=T$.
These iterations are then repeated from beginning to end until
convergence, which is reached when the computational residuals,
as defined in \cite{logg:multiadaptivity:I}, on all elements are
small enough.

\subsection{The time-slab}

The time-slab can be constructed in many ways. One is
by \emph{dyadic} partitioning, in which we compute new
time-steps for all components, based on residuals and stability
properties, choose the largest time-step $K$ as the length of the
new time-slab, and then, for every component, choose the time-step
as a fraction $K/2^n$. The advantage of such a partition are that
the time-slab has \emph{straight edges}; i.e., for every time-slab
there is a common point in time $t'$ (the end-point of the
time-slab) which is an end-point for the last element of every
component in the time-slab, and that the components have many
common nodes. The disadvantage is that the choice of time-steps
is constrained.

Another possibility is a \emph{rational} partition of the
time-slab. We choose the largest individual time-step $K$ as the
length of the time-slab, and time-steps for the remaining
components are chosen as fractions of this large time-step,
$K/2$, $K/3$, $K/4$, and so on. In this way we increase the set of
possible time-step selections, as compared to dyadic
partitioning, but the number of common nodes shared between
different components is decreased.

\begin{figure}[thbp]
\begin{center}
\leavevmode
\psfrag{1}{}
\psfrag{1/2}{}
\psfrag{1/3}{}
\psfrag{1/4}{}
\psfrag{1/5}{}
\psfrag{0}{$0$}
\psfrag{T}{$T$}
\includegraphics[width=\textwidth,height=4.5cm]{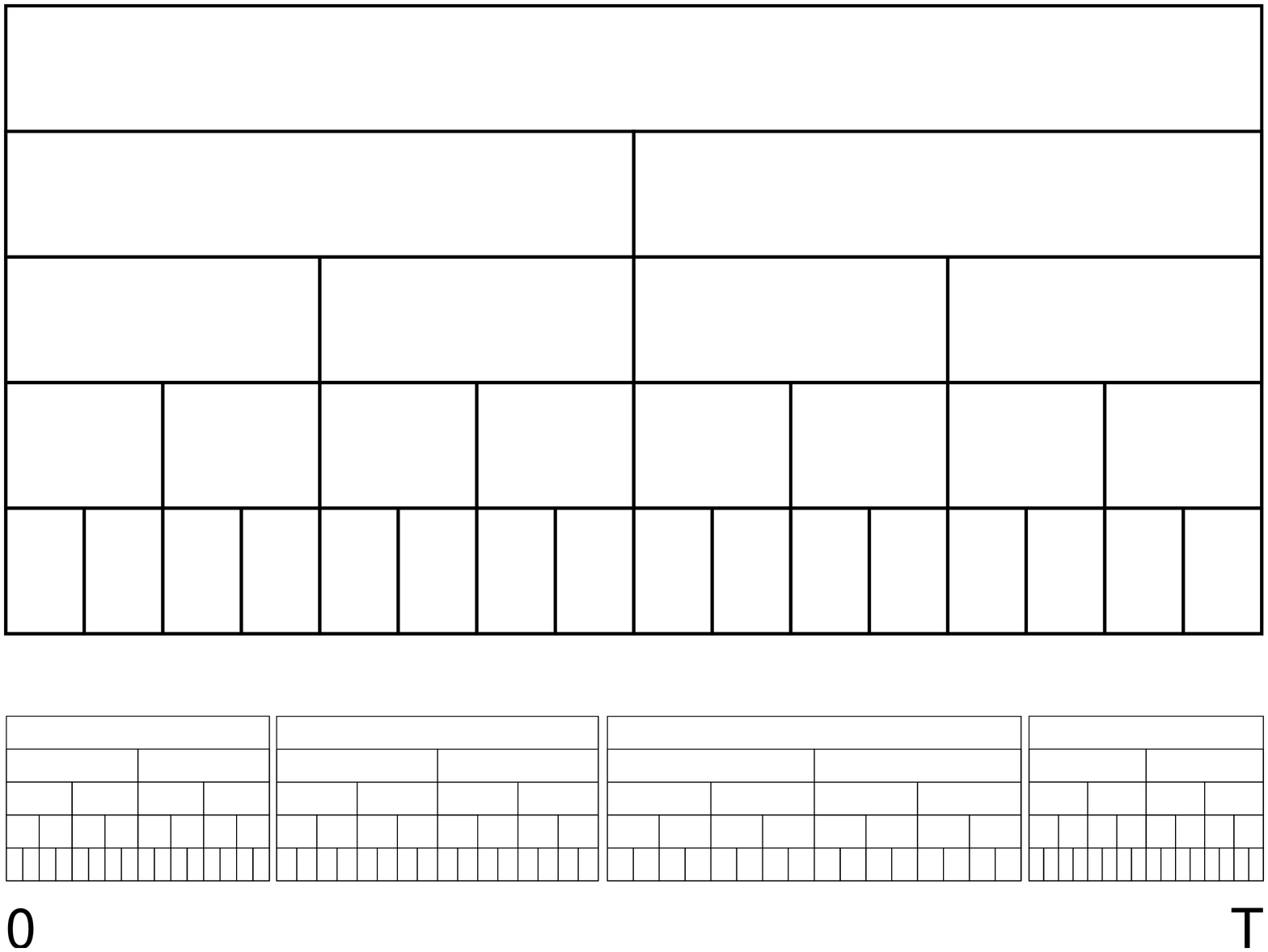}

\vspace{0.5cm}

\includegraphics[width=\textwidth,height=4.5cm]{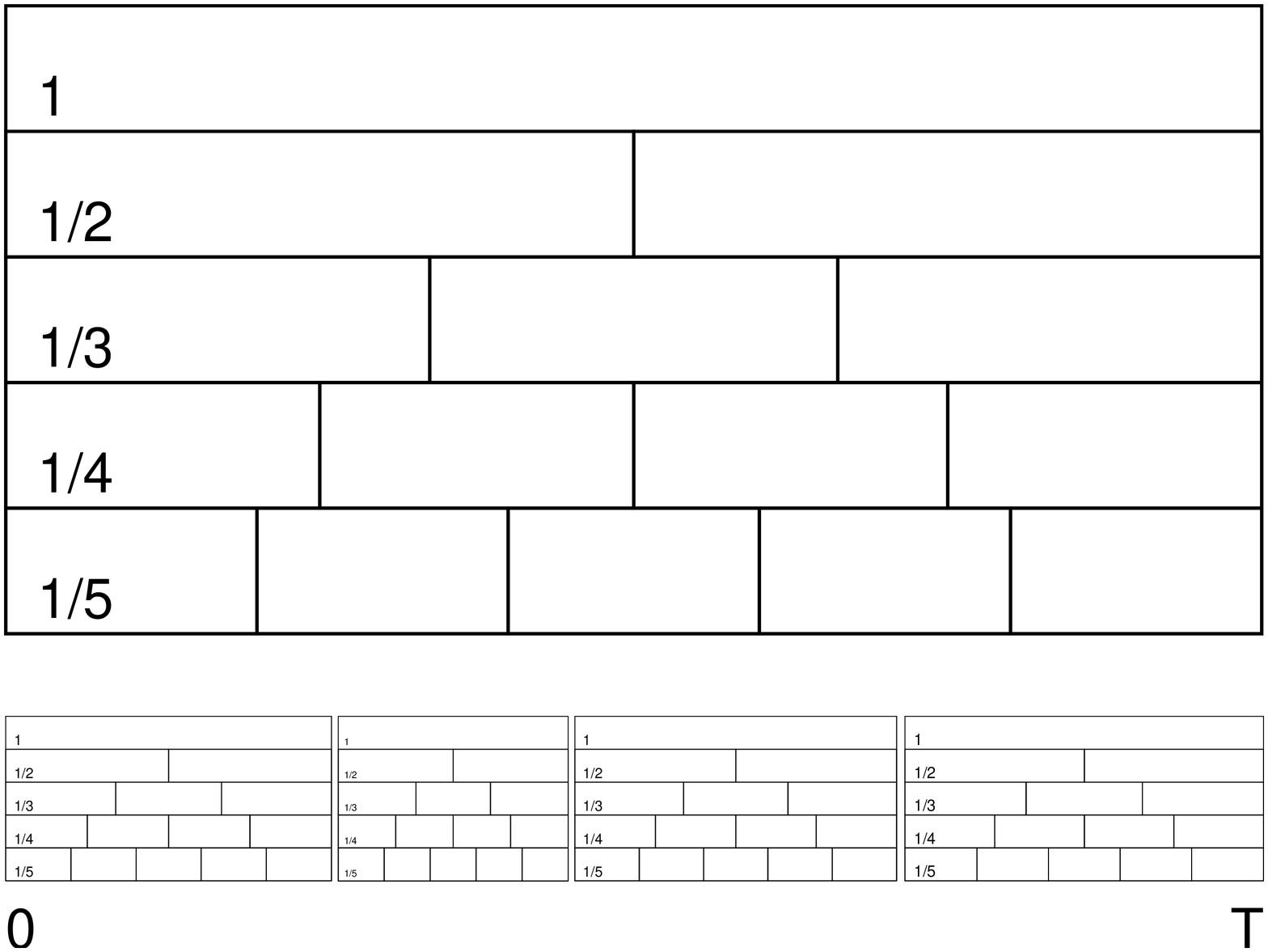}

\vspace{0.5cm}

\includegraphics[width=\textwidth,height=4.5cm]{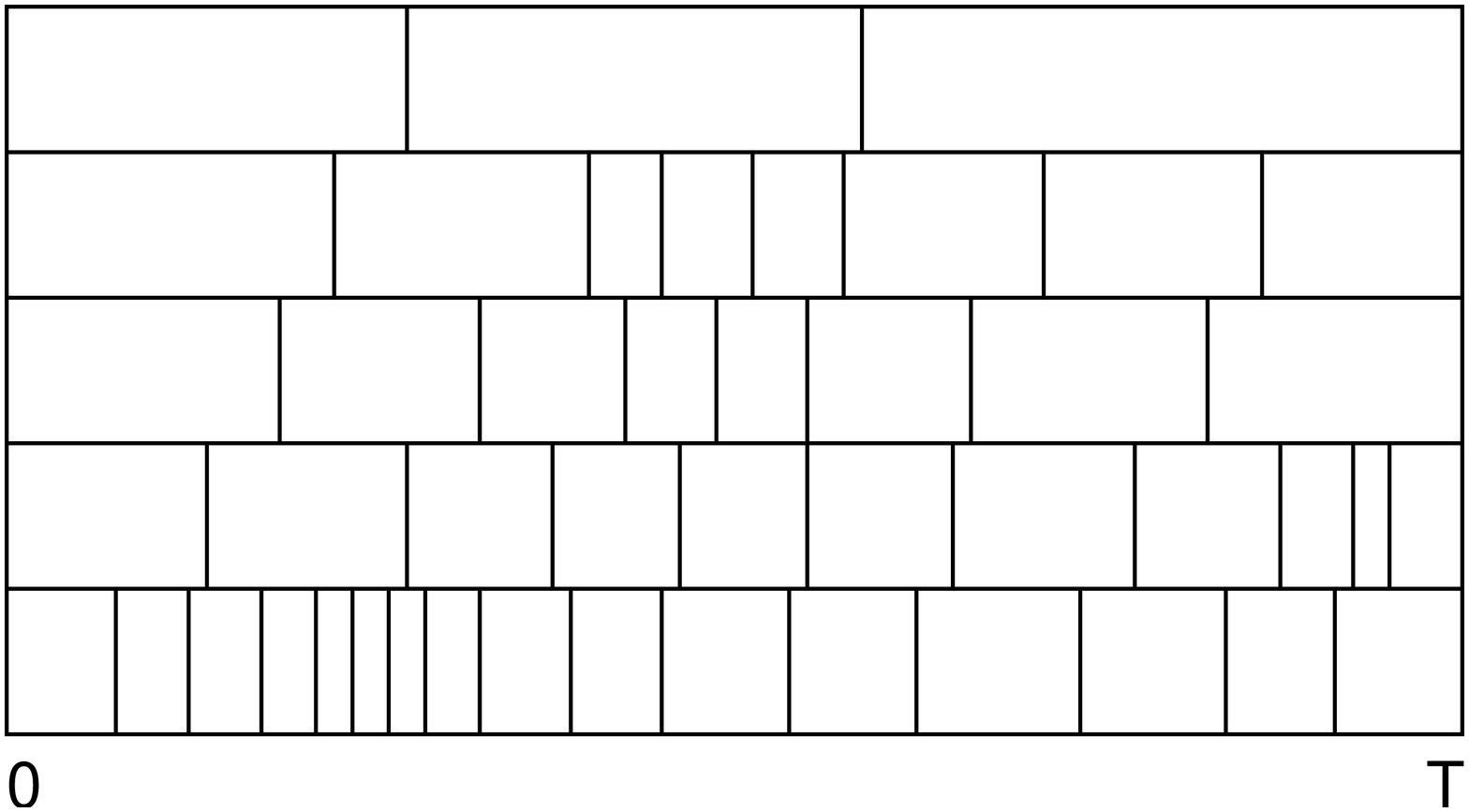}
\caption{Different choices of time-slabs. Top: a
dyadic partition of the time-slab; middle: a rational
partition; bottom: a partition used in the time-crawling
approach, where the only restriction on the time-steps is that we
match the final time end-point $T$.}
\label{fig:timeslabs}
\end{center}
\end{figure}

A third option is to not impose any constraint at all on the
selection of time-steps---except that we match the final time
end-point. The time-steps may vary also within the time-slab for
the individual components. The price we have to pay is that we
have in general no common nodes, and the edges of the time-slab
are no longer straight. We illustrate the three choices of
partitioning schemes in Figure \ref{fig:timeslabs}. Although
dyadic or rational partitioning is clearly advantageous in terms
of easier bookkeeping and common nodes, we focus below on
unconstrained time-step selection. In this way we stay close to
the original, general formulation of the multi-adaptive
methods. We refer to this as the \emph{time-crawling} approach.

\subsection{The time-crawling approach}

The construction of the time-slab brings with it a number of
technical and algorithmic problems. We will not discuss here
the implementational and data structural aspects of the
algorithm---there will be much to keep track of and this has to
be done in an efficient way---but we will give a brief account
of how the time-slab is formed and updated.

\begin{figure}[t]
\begin{center}
\leavevmode
\psfrag{0}{$0$}
\psfrag{T}{$T$}
\includegraphics[width=\textwidth,height=7cm]{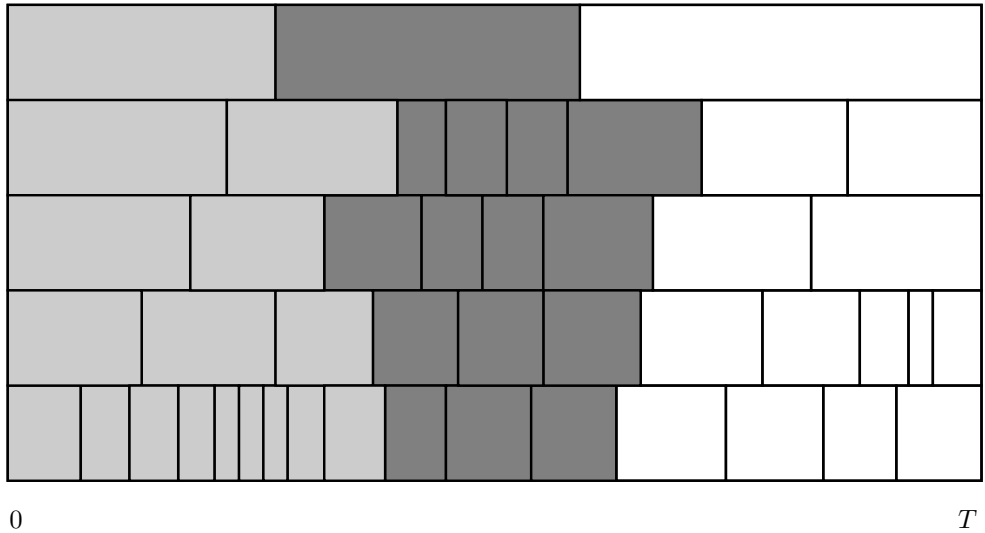}
\caption{The time-slab used in the time-crawling
approach to multi-adaptive time-stepping (dark grey). Light grey
indicates elements that have already been computed.}
\label{fig:timeslab}
\end{center}
\end{figure}

Assume that in some way we have formed a time-slab, such as the
one in Figure \ref{fig:timeslab}. We make iterations on the
time-slab, starting with the last element and continuing to the
right. After iterating through the time-slab a few times, the
computational (discrete) residuals, corresponding to the solution
of the discrete equations (\ref{eq:equations}), on all elements
have decreased below a specified tolerance for the computational
error, indicating convergence.

For the elements at the front of the slab
(those closest to time $t=T$), the values have been computed using
extrapolated values of many of the other elements. The strategy
now is to leave behind \emph{only those elements that are fully
covered by all other elements}. These are then cut off from the
time-slab, which then decreases in size. Before starting the
iterations again, we have to form a new time-slab. This will
contain the elements of the old time-slab that were not removed,
and a number of new elements. We form the new time-slab by
requiring that all elements of the previous time-slab be
totally covered within the new time-slab. In this way we know
that every new time-slab will produce at least $N$ new elements.
The time-slab is thus crawling forward in time rather than
marching.

An implementation of the method then contains the three
consecutive steps described above: iterating on an existing
time-slab, decreasing the size of the time-slab (cutting off
elements at the end of the time-slab, i.e., those closest to time $t=0$), and
incorporating new elements at the front of the time-slab.

\begin{remark}{\rm
Even if an element within the time-slab is totally covered by all
other elements, the values on this element still may not be
completely determined, if they are based on the values of some
other element that is not totally covered, or if this element is
based on yet another element that is not totally covered, and so
on. To avoid this, one can impose the requirement that the
time-slabs should have straight edges.}
\end{remark}

\subsection{Diagonal Newton}

For stiff problems the time-step condition required for
convergence of direct fixed point iteration is too restrictive,
and we need to use a more implicit solution strategy.

Applying a full Newton's method, we increase the range of allowed
time-steps and also the convergence rate, but this is costly in
terms of memory and computational time, which is especially
important in our setting, since the size of the slab may often be
much larger than the number of components, $N$ (see Figure
\ref{fig:timeslab}). We thus look for a simplified Newton's
method which does not increase the cost of solving the problem,
as compared to direct fixed point iteration, but still has some
advantages of the full Newton's method.

Consider for simplicity the case of the multi-adaptive backward
Euler method, i.e., the $\mathrm{mdG}(0)$ method with end-point
quadrature. On every element we then want to solve
\begin{equation}
U_{ij} = U_{i,j-1} + k_{ij} f_i(U(t_{ij}),t_{ij}).
\label{eq:euler}
\end{equation}
In order to apply Newton's method we write (\ref{eq:euler}) as
\begin{equation}
F(V) = 0
\label{eq:general}
\end{equation}
with $F_i(V) =   U_{ij} - U_{i,j-1} - k_{ij}
f_i(U(t_{ij}),t_{ij})$ and $V_i = U_{ij}$. Newton's method is then
\begin{equation}
V^{n+1} = V^n - (F'(V^n))^{-1} F(V^n).
\end{equation}
We now simply replace the Jacobian with its diagonal so that for
component $i$ we have
\begin{equation}
U_{ij}^{n+1} = U_{ij}^{n} - \frac{U_{ij}^n - U_{i,j-1} - k_{ij} f_i}{1-k_{ij}\PD{f_i}{u_i}}
\end{equation}
with the right-hand side evaluated at $V^n$. We now note that we
can rewrite this as
\begin{equation}
U_{ij}^{n+1} = U_{ij}^{n} - \theta (U_{ij}^n - U_{i,j-1} - k_{ij}
f_i) = (1-\theta) U_{ij}^{n} + \theta ( U_{i,j-1} +  k_{ij} f_i )
\end{equation}
with
\begin{equation}
\theta = \frac{1}{1-k_{ij}\PD{f_i}{u_i}}
\end{equation}
so that we may view the simplified Newton's method as a damped
version, with damping $\theta$, of the original fixed point
iteration.

The individual damping parameters are cheap to compute. We do not
need to store the Jacobian and we do not need linear algebra. We
still obtain some of the good properties of the full Newton's
method.

For the general \mcgq\ or \mdgq\ method, the same analysis
applies. In this case, however, when we have more degrees of
freedom to solve for on every local element,
$1-k_{ij}\PD{f_i}{u_i}$ will be a small local matrix of size
$q\times q$ for the \mcgq\ method and size $(q+1)\times (q+1)$
for the \mdgq\ method.

\subsection{Explicit or implicit}

Both \mcgq\ and \mdgq\ are implicit methods since they are
implicitly defined by the set of equations (\ref{eq:equations})
on each time-slab. However, depending on the solution strategy
for these equations, the resulting fully discrete scheme may be
of more or less explicit character. Using a diagonal Newton's
method as in the current implementation of Tanganyika, we obtain
a method of basically explicit nature. This gives an efficient
code for many applications, but we may expect to meet
difficulties for stiff problems.

\subsection{The stiffness problem}

In a stiff problem the solution varies quickly inside transients and
slowly outside transients. For accuracy the time-steps
will be adaptively kept small inside transients and then will be
within the stability limits of an explicit method, while outside
transients larger time-steps will be used.
Outside the transients the diagonal Newton's method handles stiff
problems of sufficiently diagonal nature. Otherwise the
strategy is to decrease the time-steps whenever needed for
stability reasons. Typically this results in an oscillating
sequence of time-steps where a small number of large time-steps are
followed by a small number of stabilizing small time-steps.

Our solver Tanganyika thus performs like a modern unstable jet
fighter, which needs small stabilizing wing flaps to follow a
smooth trajectory. The pertinent question is then the number of
small stabilizing time-steps per large time-step. We analyze this
question in \cite{logg:claes:stiffode} and show that for certain
classes of stiff problems it is indeed possible to successfully
use a stabilized explicit method of the form implemented in
Tanganyika.

\subsection{Preparations}

There are many ``magic numbers'' that need to be computed in
order to implement the multi-adaptive methods, such as quadrature
points and weights, the polynomial weight functions evaluated at
these quadrature points, etc. In Tanganyika, these numbers are
computed at the startup of the program and stored for efficient
access. Although somewhat messy to compute, these are all
computable by standard techniques in numerical analysis; see, e.g.,
\cite{book:numericalrecipes}.

\subsection{Solving the dual problem}

In addition to solving the primal problem (\ref{eq:u'=f}), we
also have to solve the continuous dual problem to obtain error
control. This is an ODE in itself that we can solve using the
same solver as for the primal problem.

In order to solve this ODE, we need to know the Jacobian of $f$
evaluated at a mean value of the true solution $u(t)$ and the
approximate solution $U(t)$. If $U(t)$ is sufficiently close to
$u$, which we will assume, we approximate the (unknown) mean
value by $U(t)$. When solving the dual, the primal solution must
be accessible, and the Jacobian must be computed numerically by
difference quotients if it is not explicitly known. This makes
the computation of the dual solution expensive. Error control can,
however, be obtained at a reasonable cost: for one thing, the dual
problem does not have to be solved with as high a precision as the
primal problem; a relative error of, say, $10\%$ may be disastrous
for the primal, whereas for the dual this only means that the
error estimate will be off by $10\%$, which is acceptable.
Second, the dual problem is linear, which may be taken into
account when implementing a solver for the dual. If we can afford
the linear algebra, as we can for reasonably small systems, we
can solve the discrete equations directly without any iterations.


\section{Tanganyika}

We now give a short description of the implementation of the
multi-adaptive methods, \emph{Tanganyika}, which has been used to
obtain the numerical results presented below.

\subsection{Purpose}

The purpose of Tanganyika \cite{logg:tanganyika} is to be a
working implementation of the multi-adaptive methods.
The code is open-source (GNU GPL \cite{www:gnu}), which
means that anyone can freely review the code, which is
available at
http://www.phi. chalmers.se/tanganyika/.
Comments are welcome.

\subsection{Structure and implementation}

The solver is implemented as a C/C++ library.
The C++ language makes abstraction easy, which allows the
implementation to follow closely the formulation of the two
methods. Different objects and algorithms are thus implemented as
C++ classes, including
\texttt{Solution}, \texttt{Element}, \texttt{cGqElement}, \texttt{dGqElement},
\texttt{TimeSlab},
\texttt{ErrorControl}, \texttt{Galerkin}, \texttt{Component},
and so on.

\section{Applications}

In this section, we present numerical results for a variety of
applications. We discuss some of the problems in detail
and give a brief overview of the rest. A more extensive account
can be found in \cite{logg:lic:II}.

\subsection{A simple test problem}
\label{sec:testproblem}

\begin{figure}[t]
\begin{center}
\includegraphics[width=\textwidth]{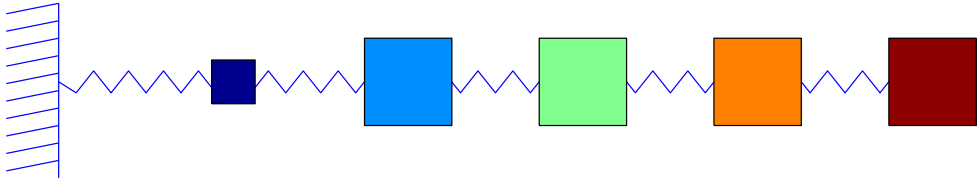}
\caption{A mechanical system consisting of $N=5$ masses attached with springs.}
\label{fig:multiadaptivity-system}
\end{center}
\end{figure}

To demonstrate the potential of the multi-adaptive methods, we
consider a dynamical system in which a small part of the system
oscillates rapidly. The problem is to accurately compute the
positions (and velocities) of the $N$ point-masses attached
with springs of equal stiffness, as in Figure
\ref{fig:multiadaptivity-system}.

If we choose one of the masses to be much smaller than the others, $m_1
= 10^{-4}$ and $m_2=m_3=\cdots =m_N=1$, then we expect the dynamics of
the system to be dominated by the smallest mass, in the sense
that the resolution needed to compute the solution will be
completely determined by the fast oscillations of the smallest
mass.

To compare the multi-adaptive method with a standard method, we
first compute with constant time-steps $k=k_0$ using the standard
$\mathrm{cG}(1)$ method and measure the error, the cpu time
needed to obtain the solution, the total number of steps, i.e., $M
= \sum_{i=1}^N M_i$, and the number of local function evaluations.
We then compute the solution with individual time-steps, using
the $\mathrm{mcG}(1)$ method, choosing the time-steps $k_i =
k_0$ for the position and velocity components of the smallest
mass, and choosing $k_i = 100 k_0$ for the other components (knowing
that the frequency of the oscillations scales like $1/\sqrt{m}$).
For demonstration purposes, we thus choose the time-steps a
priori to fit the dynamics of the system.

\begin{figure}[thbp]
\begin{center}
\psfrag{e}{\hspace{-0.7cm}$\|e(T)\|_{\infty}$}
\psfrag{cpu time}{\hspace{-1.0cm}\textsf{cpu time / seconds}}
\psfrag{function evaluations}{\hspace{-0.5cm}\textsf{function evaluations}}
\psfrag{steps}{\textsf{steps}}
\includegraphics[width=\textwidth,height=9.8cm]{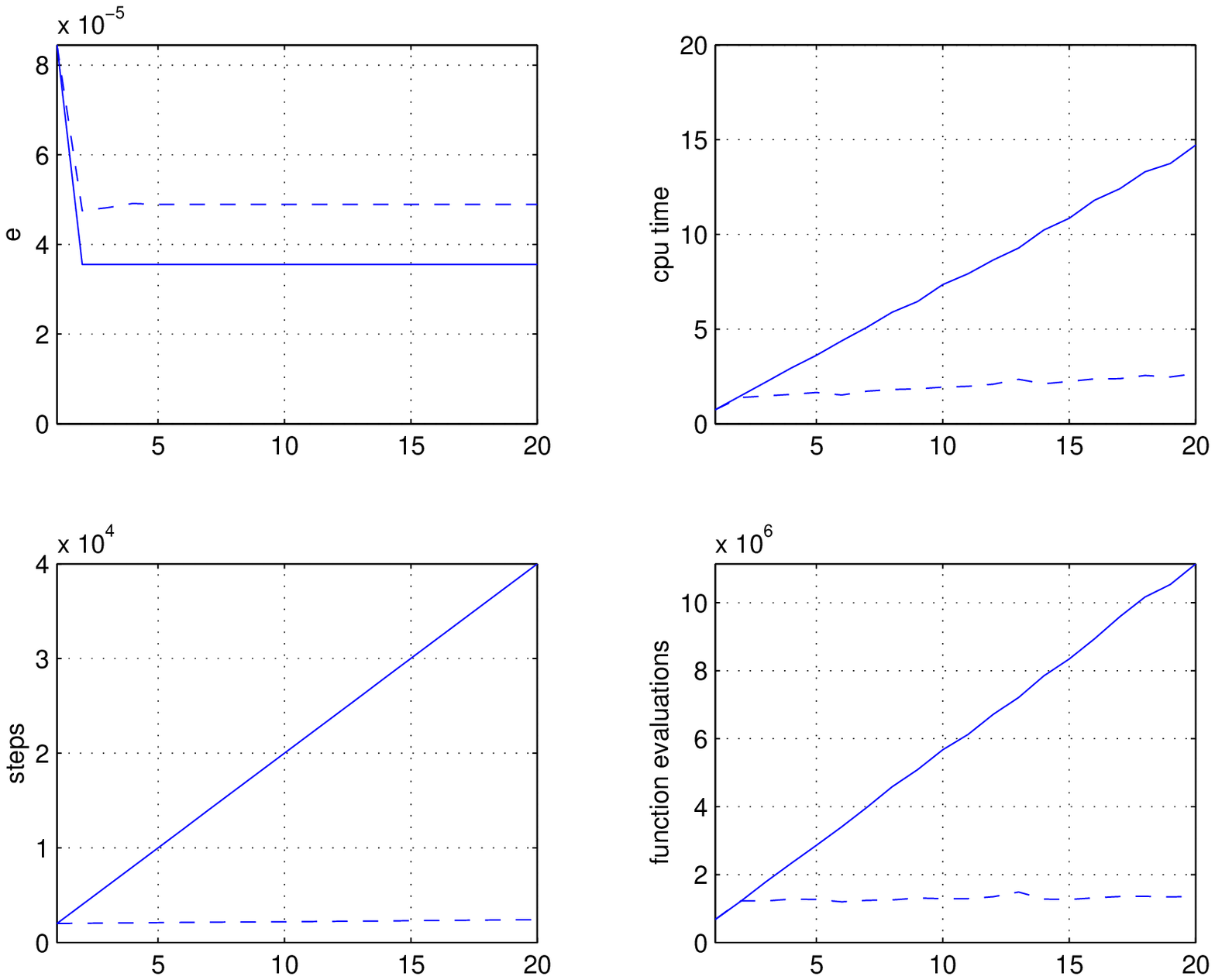}
\includegraphics[width=\textwidth,height=9.8cm]{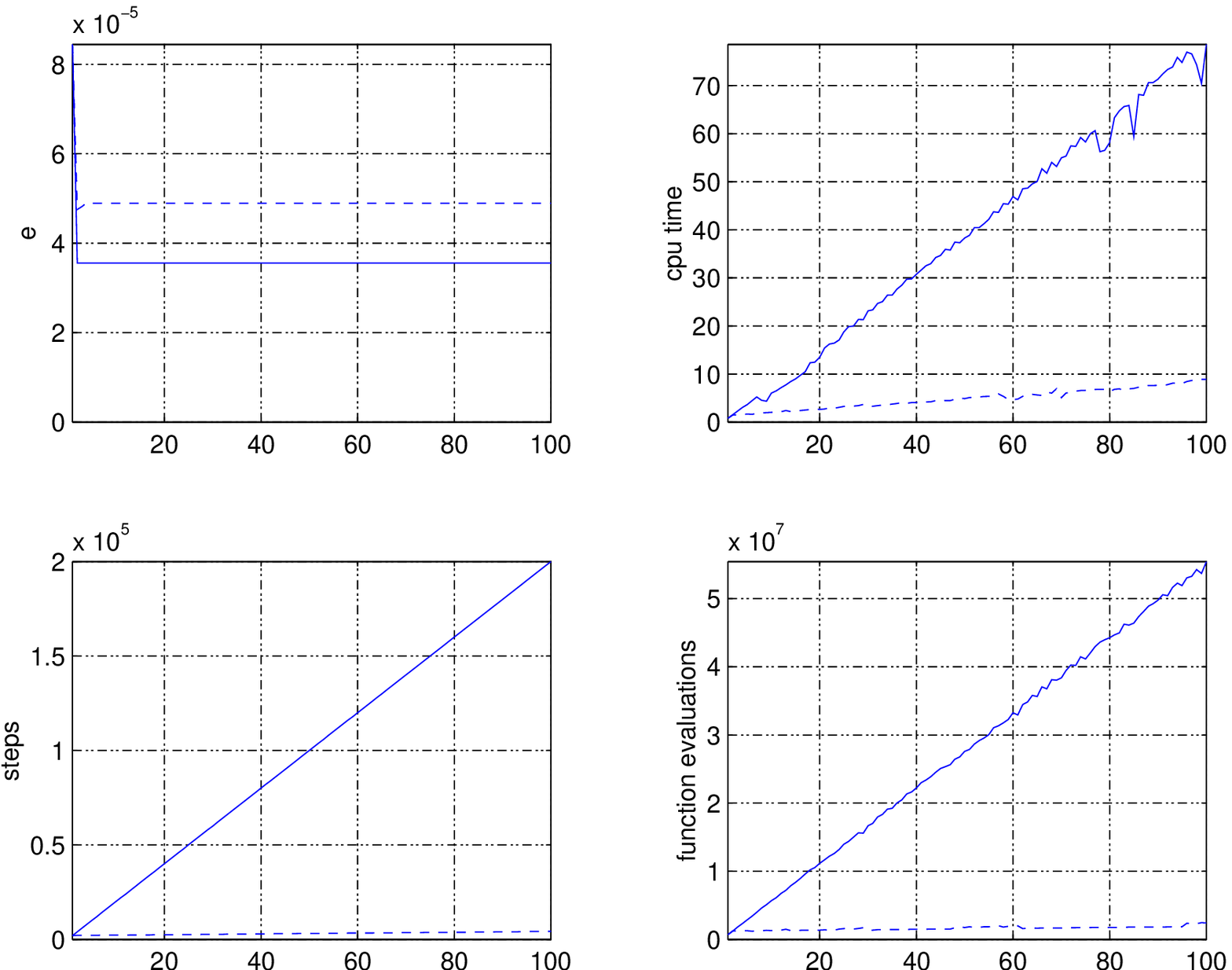}
\caption{Error, cpu time, total number of steps, and
number of function evaluations as function of the number of
masses for the multi-adaptive $\mathrm{cG}(1)$ method (dashed lines)
and the standard $\mathrm{cG}(1)$ method (solid lines).}
\label{fig:multiadaptivity-results}
\end{center}
\end{figure}

We repeat the experiment for increasing values of $N$ (see Figure
\ref{fig:multiadaptivity-results}) and find that
the error is about the same and constant for both methods.  As $N$
increases, the total number of time-steps, the number of local
function evaluations (including also residual evaluations), and
the cpu time increase linearly for the standard method, as we
expect.  For the multi-adaptive method, on the other hand, the
total number of time-steps and local function evaluations remains
practically constant as we increase $N$. The cpu time increases
somewhat, since the increasing size of the time-slabs introduces
some overhead, although not nearly as much as for the standard
method.  For this particular problem the gain of the
multi-adaptive method is thus a factor $N$, where $N$ is the size
of the system, so that by considering a large-enough system, the
gain is arbitrarily large.

\subsection{The Lorenz system}
\label{sec:lorenz}

We consider now the famous Lorenz system,
\begin{equation}
\left\{
\begin{array}{rcl}
\dot{x} &=& \sigma (y-x), \\
\dot{y} &=& rx - y -xz,   \\
\dot{z} &=& xy - bz,      \\
\end{array}
\right.
\label{eq:lorenz}
\end{equation}
with the usual data $(x(0),y(0),z(0))=(1,0,0)$, $\sigma=10$,
$b=8/3$, and $r=28$; see \cite{EstJoh98}. The solution
$u(t)=(x(t),y(t),z(t))$ is very sensitive to perturbations and is
often described as ``chaotic.'' With our perspective this
is reflected by stability factors with rapid growth in time.\pagebreak

\begin{figure}[t]
\begin{center}
\leavevmode
\psfrag{x}{$x$}
\psfrag{z}{$z$}
\psfrag{U}{$U$}
\psfrag{t}{$t$}
\myfigure{width=\textwidth}{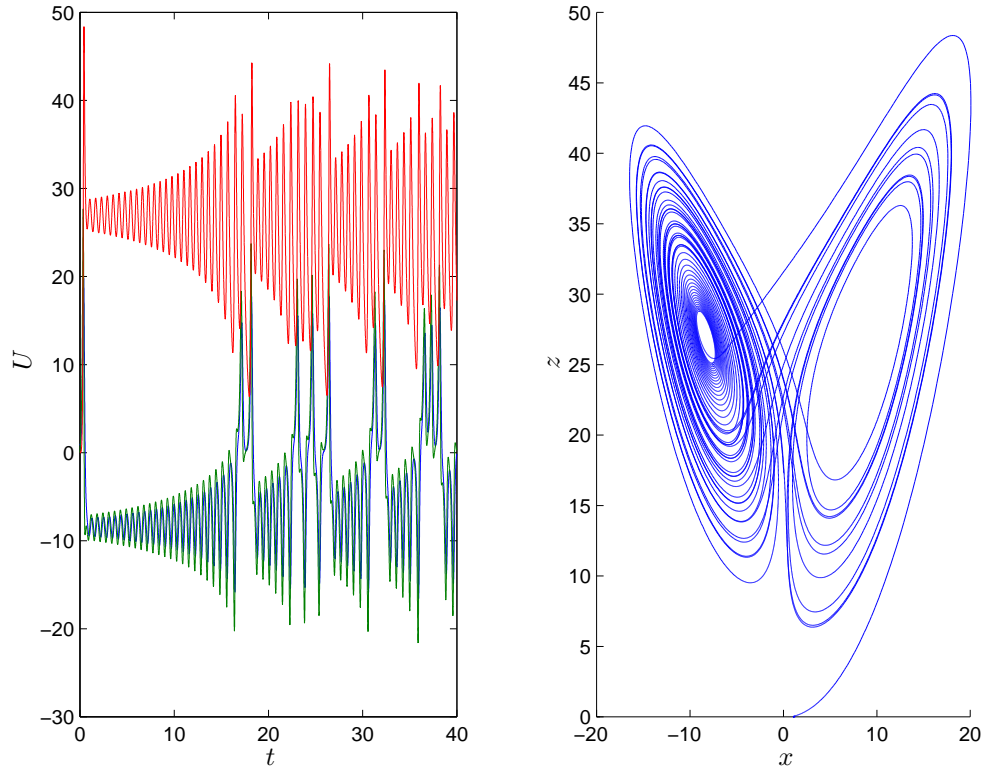}
\caption{On the right is the trajectory of the Lorenz system for
final time $T=40$, computed with the multi-adaptive
$\mathrm{cG}(5)$ method. On the left is a plot of the time
variation of the three components.}
\label{fig:lorenz-trajectory}
\end{center}
\end{figure}

The computational challenge is to solve the Lorenz system
accurately on a time-interval $[0,T]$ with $T$ as large as
possible. In Figure \ref{fig:lorenz-trajectory} is shown a computed solution which is accurate on the interval $[0,40]$.
We investigate the computability of the Lorenz system
by solving the dual problem and computing stability factors to
find the maximum value of $T$. The focus in this section is
not on multi-adaptivity---we will use the same time-steps for
all components, and so
\mcgq\ becomes $\mathrm{cG}(q)$---but on higher-order methods
and the precise information that can
be obtained about the computability of a system from solving the
dual problem and computing stability factors.

\begin{figure}[thbp]
\begin{center}
\leavevmode
\psfrag{cG(5)}{\hspace{-0.15cm}$\mathrm{cG}(5)$}
\psfrag{cG(6)}{\hspace{-0.15cm}$\mathrm{cG}(6)$}
\psfrag{cG(7)}{\hspace{-0.15cm}$\mathrm{cG}(7)$}
\psfrag{cG(8)}{\hspace{-0.15cm}$\mathrm{cG}(8)$}
\psfrag{cG(9)}{\hspace{-0.15cm}$\mathrm{cG}(9)$}
\psfrag{cG(10)}{\hspace{-0.15cm}$\mathrm{cG}(10)$}
\psfrag{cG(11)}{\hspace{-0.15cm}$\mathrm{cG}(11)$}
\psfrag{cG(12)}{\hspace{-0.15cm}$\mathrm{cG}(12)$}
\psfrag{cG(13)}{\hspace{-0.15cm}$\mathrm{cG}(13)$}
\psfrag{cG(14)}{\hspace{-0.15cm}$\mathrm{cG}(14)$}
\psfrag{cG(15)}{\hspace{-0.15cm}$\mathrm{cG}(15)$}
\psfrag{t}{}
\myfigure{width=12.5cm}{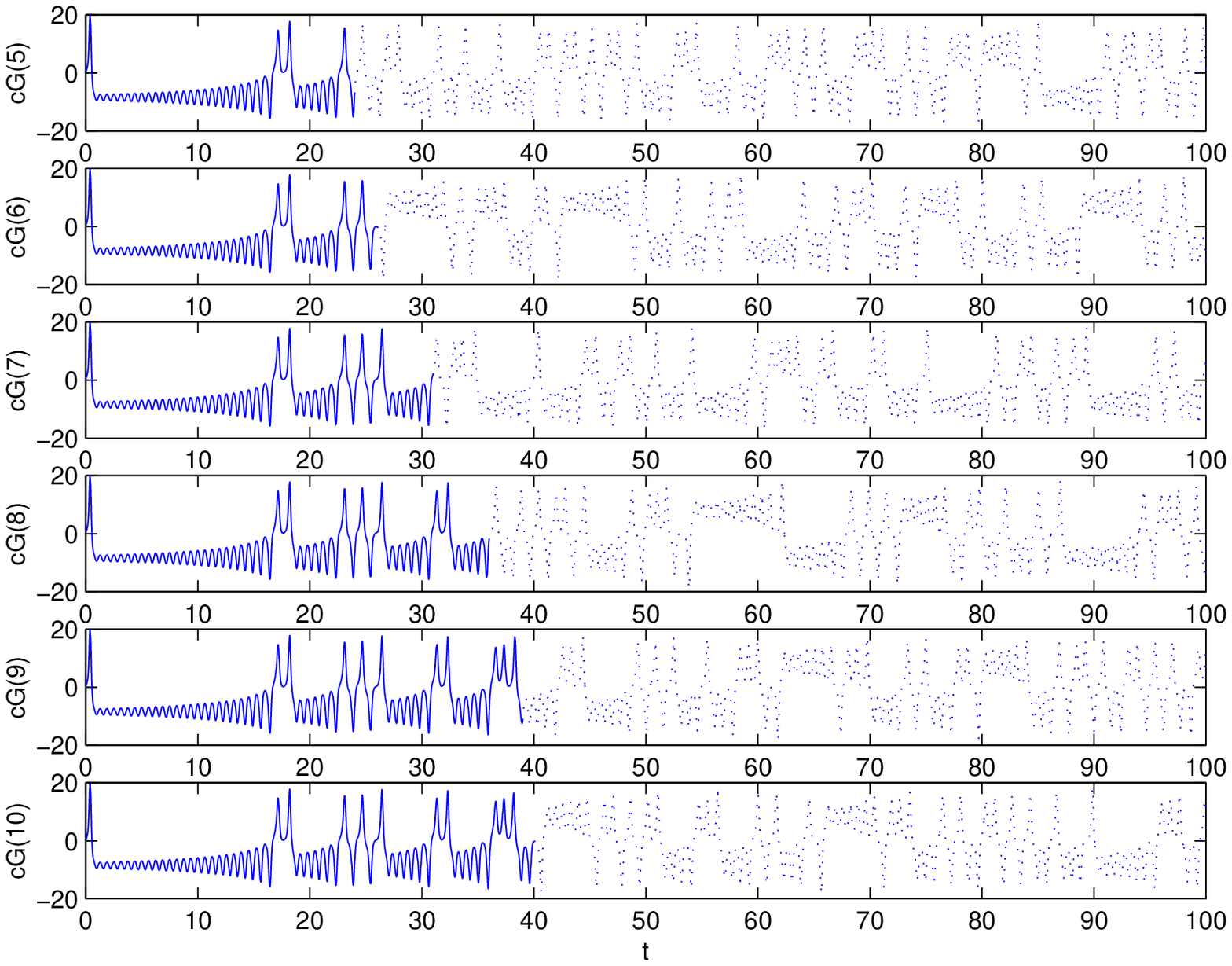}\\
\myfigure{width=12.5cm}{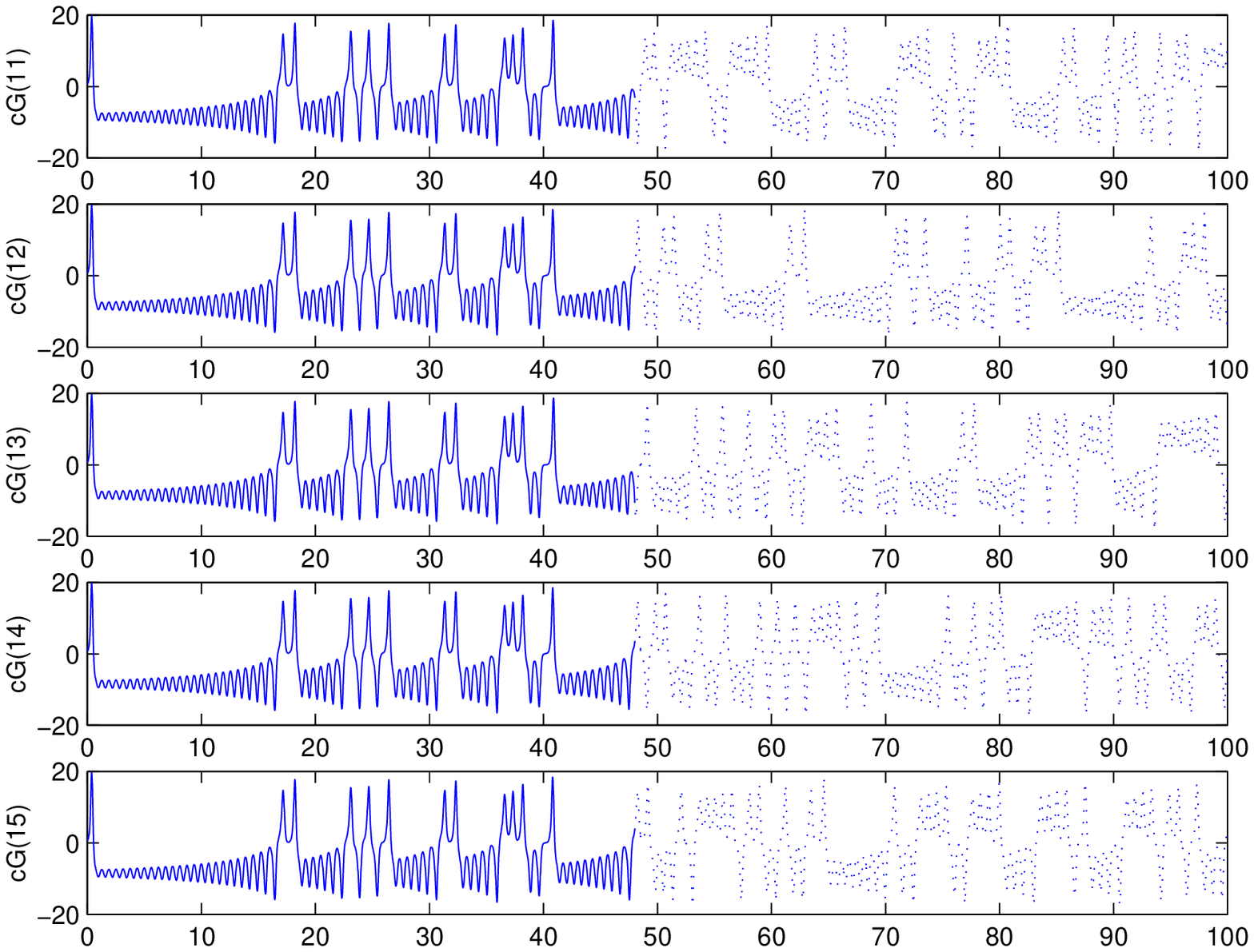}
\caption{Solutions for the $x$-component of the Lorenz system
with methods of different order, using a constant time-step
$k=0.1$. Dotted lines indicate the point beyond which the
solution is no longer accurate.}
\label{fig:lorenz-higher_order}
\end{center}
\end{figure}

As an illustration, we present in Figure
\ref{fig:lorenz-higher_order} solutions obtained with different
methods and constant time-step $k = 0.1$ for all components. For
the lower-order methods, $\mathrm{cG}(5)$ to $\mathrm{cG}(11)$,
it is clear that the error decreases as we increase the order.
Starting with the $\mathrm{cG}(12)$ method, however, the error
does not decrease as we increase the order. To explain this, we
note that in every time-step a small round-off error of size
$\sim 10^{-16}$ is introduced if the computation is carried out
in double precision arithmetic. These errors accumulate at a rate
determined by the growth of the stability factor for the
computational error (see \cite{logg:multiadaptivity:I}). As we
shall see below, this stability factor grows exponentially for
the Lorenz system and reaches a value of $10^{16}$ at final time
$T=50$, and so at this point the accumulation of round-off errors
results in a large computational error.

\subsubsection{Stability factors}

\begin{figure}[thbp]
\begin{center}
\leavevmode
\psfrag{t}{$T$}
\psfrag{S0}{$\bar{S}^{[0]}(T)$}
\psfrag{25}{}
\myfigure{width=\textwidth}{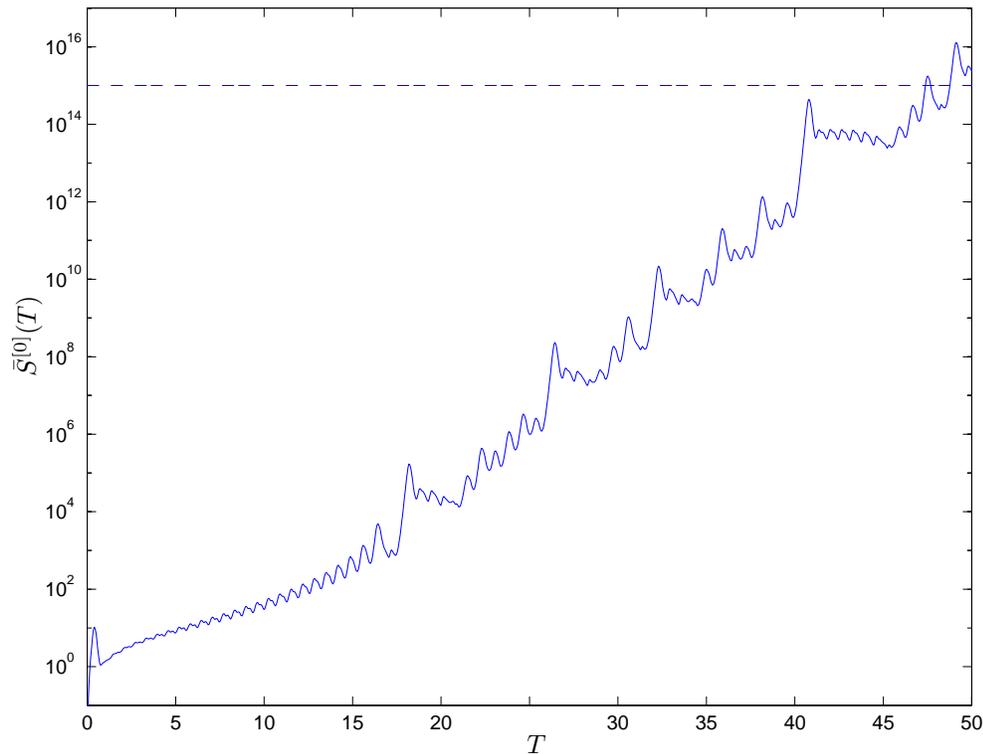}
\caption{The stability factor for computational and
quadrature errors, as function of time for the Lorenz system.}
\label{fig:lorenz-stabilityfactor_0}
\end{center}
\end{figure}
We now investigate the computation of stability factors for the
Lorenz system. For simplicity we focus on global stability
factors, such as
\begin{equation}
S^{[q]}(T) = \max_{\|v\|=1} \int_0^T \| \varphi^{(q)}(t) \| \ dt,
\label{eq:stabilityfactors}
\end{equation}
where $\varphi(t)$ is the solution of the dual problem obtained
with $\varphi(T)=v$ (and $g=0$). Letting $\Phi(t)$ be the
fundamental solution of the dual problem, we have
\begin{equation}
\  \ \,\qquad\max_{\|v\|=1} \int_0^T \| \Phi^{(q)}(t) v \| \ dt \leq
\int_0^T \max_{\|v\|=1} \| \Phi^{(q)}(t) v \| \ dt =
\int_0^T \| \Phi^{(q)}(t) \| \ dt\! =\! \bar{S}^{[q]}(T).
\label{eq:stabilityfactors,estimate}
\end{equation}
This gives a bound $\bar{S}^{[q]}(T)$ for $S^{[q]}(T)$, which for
the Lorenz system turns out to be quite sharp and which is
simpler to compute since we do not have to compute the maximum in
(\ref{eq:stabilityfactors}).

In Figure \ref{fig:lorenz-stabilityfactor_0} we plot the growth
of the stability factor for $q=0$, corresponding to computational
and quadrature errors as described in
\cite{logg:multiadaptivity:I}. The stability factor grows
exponentially with time, but not as fast as indicated by an a
priori error estimate. An a priori error estimate indicates that
the stability factors grow as
\begin{equation}
S^{[q]}(T) \sim \mathcal{A}^q e^{\mathcal{A} T},
\end{equation}
where $\mathcal{A}$ is a bound for the Jacobian of the right-hand
side for the Lorenz system. A simple estimate is
$\mathcal{A}=50$, which already at $T=1$ gives $S^{[0]}(T)\approx
10^{22}$. In view of this, we would not be able to compute even
to $T=1$, and certainly not to $T=50$, where we have
$S^{[0]}(T)\approx 10^{1000}$. The point is that although the
stability factors grow very rapidly on some occasions, such as
near the first flip at $T=18$, the growth is not monotonic. The
stability factors thus \emph{effectively} grow at a moderate
exponential~rate.

\subsubsection{Conclusions}

To predict the computability of the Lorenz system, we estimate
the growth rate of the stability factors. A simple approximation
of this growth rate, obtained from numerically computed solutions
of the dual problem, is
\begin{equation}
\bar{S}^{[q]}(T) \approx 4 \cdot 10^{(q-3)+0.37T}
\label{eq:approximation,better}
\end{equation}
or just
\begin{equation}
\bar{S}^{[q]}(T) \approx 10^{q+T/3}.
\label{eq:approximation,simpler}
\end{equation}
From the a posteriori error estimates presented in \cite{logg:multiadaptivity:I}, we find that
the computational error can be estimated as
\begin{equation}
E_C \approx S^{[0]}(T) \max_{[0,T]} \|\mathcal{R}^{\mathcal{C}}\|,
\end{equation}
where the computational residual $\mathcal{R}^{\mathcal{C}}$ is
defined as
\begin{equation}
\mathcal{R}^{\mathcal{C}}_{i}(t) = \frac{1}{k_{ij}} \left( U(t_{ij}) - U(t_{i,j-1}) - \int_{I_{ij}} f_i(U,\cdot) \ dt \right).
\end{equation}
With 16 digits of precision a simple estimate for the
computational residual is $\frac{1}{k_{ij}}\! 10^{-16}$, which
gives the approximation
\begin{equation}
\label{eq:estimate1}
E_C \approx 10^{T/3} \frac{1}{\min k_{ij}}
10^{-16} = 10^{T/3-16} \frac{1}{\min k_{ij}}.
\end{equation}
With time-steps $k_{ij}=0.1$ as above we then have
\begin{equation}
\label{eq:estimate}
E_C \approx 10^{T/3-15},
\end{equation}
and so already at time $T=45$ we have $E_C\approx 1$ and the
solution is no longer accurate.
We thus conclude by examination of the stability factors that it
is difficult to reach beyond time $T=50$ in double precision
arithmetic. (With quadruple precision we would be able to reach
time $T=100$.)

\subsection{The solar system}
\label{sec:solar}

We now consider the solar system, including the Sun, the Moon,
and the nine planets, which is a particularly important $n$-body
problem of the form
\begin{equation}
m_i \ddot{x}_{i} = \sum_{j\neq i} \frac{G m_i m_j}{|x_j - x_i|^3} (x_j - x_i),
\end{equation}
where $x_i(t)=(x_i^1(t),x_i^2(t),x_i^3(t))$ denotes the position of body $i$ at time $t$,
$m_i$ is the mass of body $i$, and $G$ is the gravitational constant.

As initial conditions we take the values at 00.00 Greenwich mean time
on January 1, 2000, obtained from the United States Naval Observatory
\cite{www:usnavy}, with initial velocities obtained by fitting a
high-degree polynomial to the values of December 1999. This
initial data should be correct to five or more digits, which is
similar to the available precision for the masses of the planets.
We normalize length and time to have the space coordinates per
astronomical unit, AU, which is (approximately) the mean distance
between the Sun and Earth, the time coordinates per year, and the
masses per solar mass. With this normalization, the gravitational
constant is $4\pi^2$.

\subsubsection{Predictability}

\begin{figure}[t]
\begin{center}
\psfrag{v}{$\varphi$}
\psfrag{t}{$t$}
\includegraphics[width=\textwidth]{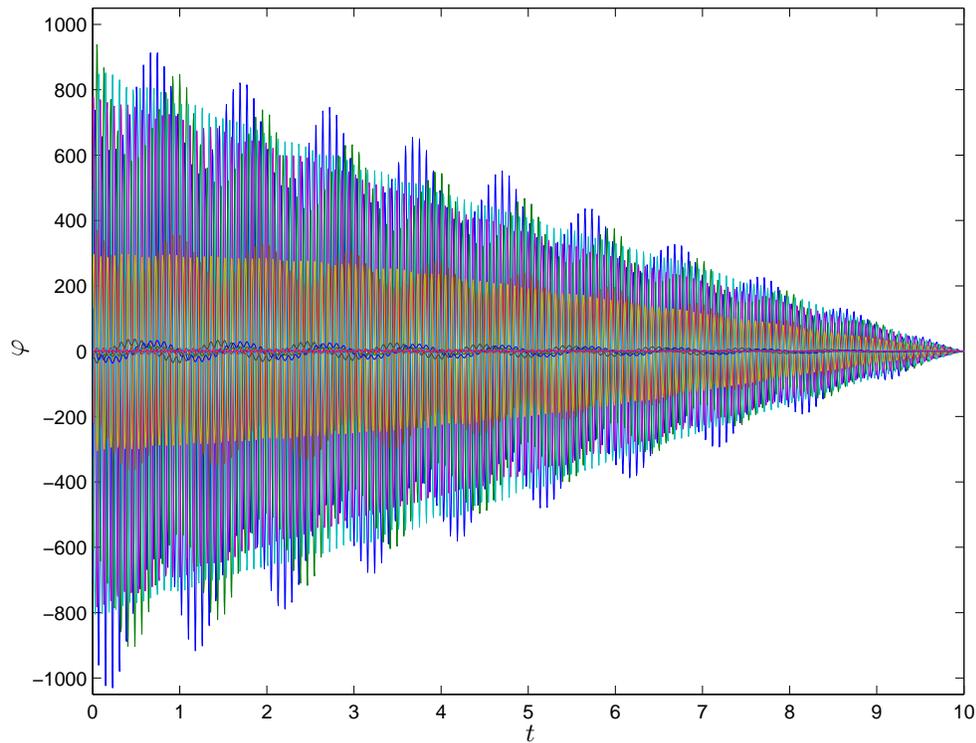}
\caption{Part of the dual of the solar system with
data chosen for control of the error in the position of the moon
at final time.}
\label{fig:earth-moon-dual}
\end{center}
\end{figure}

Investigating the \emph{predictability} of the solar system, the
question is how far we can accurately compute the solution, given
the precision in initial data. In order to predict the
accumulation rate of errors, we solve the dual problem and
compute stability factors. Assuming the initial data is correct
to five or more digits, we find that the solar system is
computable on the order of 500 years. Including also the Moon, we
cannot compute more than a few years. The dual solution grows
linearly backward in time (see Figure \ref{fig:earth-moon-dual}),
and so errors in initial data grow linearly with time. This means
that for every extra digit of increased precision, we can reach
10 times further.

\subsubsection{Computability}

To touch briefly on the fundamental question of the
\emph{computability} of the solar system, concerning how far the
system is computable with correct initial data and correct model,
we compute the trajectories for Earth, the Moon, and the Sun over
a period of $50$ years, comparing different methods. Since errors
in initial data grow linearly, we expect numerical errors as well
as stability factors to grow quadratically.

\begin{figure}[thbp]
\begin{center}
\psfrag{cG1}{$\mathrm{cG}(1)$}
\psfrag{cG2}{$\mathrm{cG}(2)$}
\psfrag{dG1}{$\mathrm{dG}(1)$}
\psfrag{dG2}{$\mathrm{dG}(2)$}
\psfrag{e}{$e$}
\psfrag{t}{$t$}
\includegraphics[width=\textwidth]{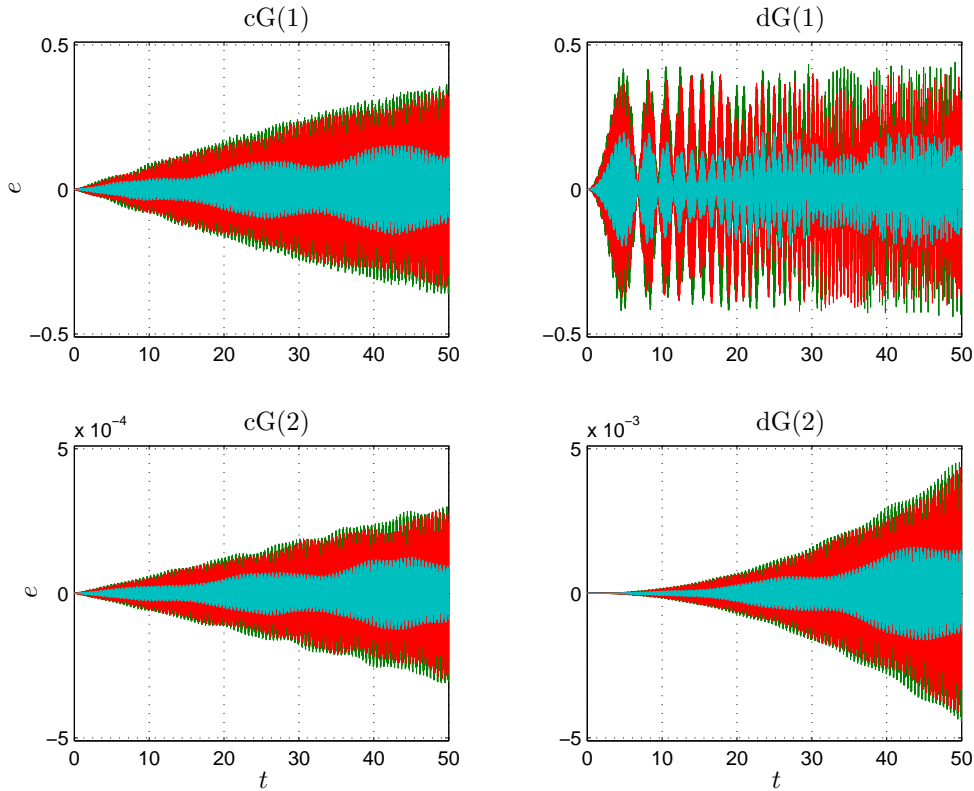}
\caption{The growth of the error over $50$ years for the
Earth--Moon--Sun system as described in the text.}
\label{fig:earth-moon-computability}
\end{center}
\end{figure}

In Figure \ref{fig:earth-moon-computability} we plot the errors
for the $18$ components of the solution, computed for $k=0.001$
with $\mathrm{cG}(1)$, $\mathrm{cG}(2)$, $\mathrm{dG}(1)$, and
$\mathrm{dG}(2)$. This figure contains much information. To
begin with, we see that the error seems to grow linearly for the
$\mathrm{cG}$ methods. This is in accordance with earlier
observations \cite{estep:lineargrowth,mats:lineargrowth} for
periodic Hamiltonian systems, recalling that the
$\mathrm{(m)cG}(q)$ methods conserve energy
\cite{logg:multiadaptivity:I}. The stability factors, however,
grow quadratically and thus overestimate the error growth for
this particular problem. In an attempt to give an intuitive
explanation of the linear growth, we may think of the error
introduced at every time-step by an energy-conserving method as a
pure phase error, and so at every time-step the Moon is pushed
slightly forward along its trajectory (with the velocity adjusted
accordingly). Since a pure phase error does not accumulate but
stays constant (for a circular orbit), the many small phase
errors give a total error that grows linearly with time.

Examining the solutions obtained with the $\mathrm{dG}(1)$ and
$\mathrm{dG}(2)$ methods, we see that the error grows
quadratically, as we expect. For the $\mathrm{dG}(1)$ solution,
the error reaches a maximum level of $\sim 0.5$ for the velocity
components of the Moon. The error in position for the Moon is
much smaller. This means that the Moon is still in orbit around
Earth, the position of which is still very accurate, but the
position relative to Earth, and thus also the velocity,
is incorrect. The error thus grows quadratically until it reaches a
limit. This effect is also visible for the error of the
$\mathrm{cG}(1)$ solution; the linear growth flattens out as the
error reaches the limit. Notice also that even if the
higher-order $\mathrm{dG}(2)$ method performs better than the
$\mathrm{cG}(1)$ method on a short time-interval, it will be
outrun on a long enough interval by the $\mathrm{cG}(1)$ method,
which has linear accumulation of errors (for this particular
problem).

\begin{figure}[t]
\begin{center}
\psfrag{e}{$e$}
\psfrag{k}{$k$}
\psfrag{t}{$t$}
\psfrag{t1}{\footnotesize All components}
\psfrag{t2}{\footnotesize Position of the Moon}
\includegraphics[width=\textwidth]{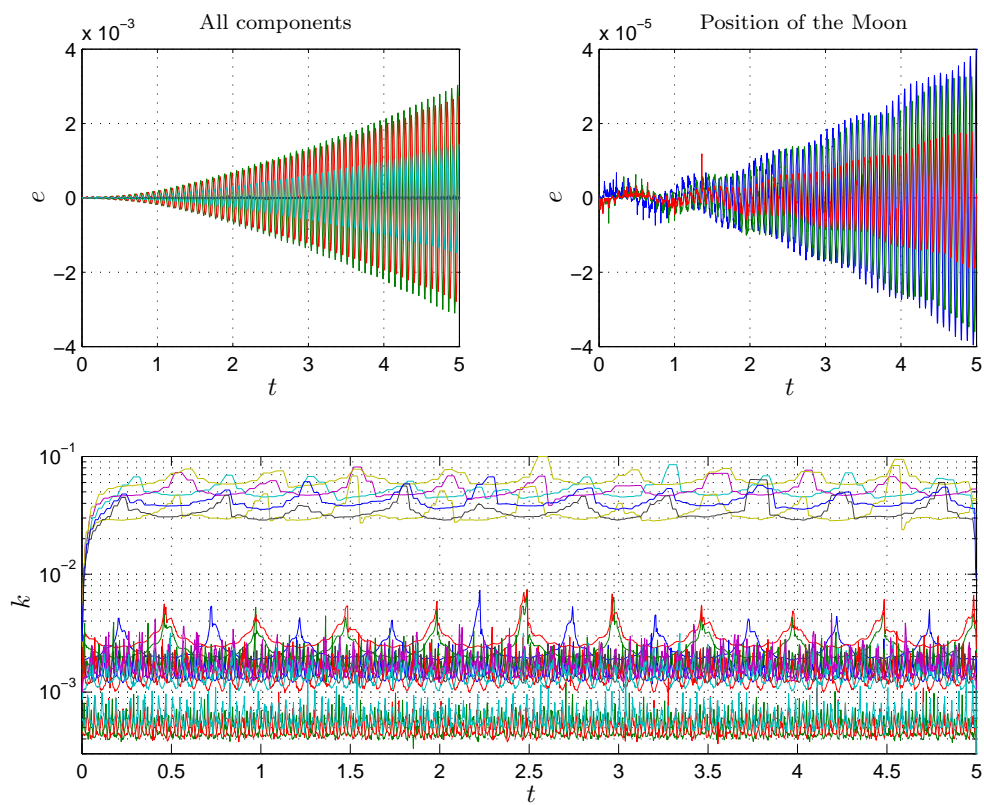}
\caption{The growth of the error over $5$ years for the
Earth--Moon--Sun system computed with the $\mathrm{mcG}(2)$ method,
together with the multi-adaptive time-steps.}
\label{fig:earth-moon-computability-multiadaptive}
\end{center}
\end{figure}

\subsubsection{Multi-adaptive time-steps}

Solving with the multi-adaptive method $\mathrm{mcG}(2)$ (see
Figure \ref{fig:earth-moon-computability-multiadaptive}), the
error grows quadratically. We saw in
\cite{logg:multiadaptivity:I} that in order for the
$\mathrm{mcG}(q)$ method to conserve energy, we require that
corresponding position and velocity components use the same
time-steps. Computing with different time-steps for all
components, as here, we thus cannot expect to have linear error
growth. Keeping $k_i^2 r_i \leq \mathrm{tol}$ with $\mathrm{tol}
= 10^{-10}$ as here, the error grows as $10^{-4} T^2$ and we are
able to reach $T\sim 100$. Decreasing $\mathrm{tol}$ to, say,
$10^{-18}$, we could instead reach $T\sim 10^6$.
\begin{figure}[thbp]
\begin{center}
\leavevmode
\includegraphics[width=\textwidth,height=10cm]{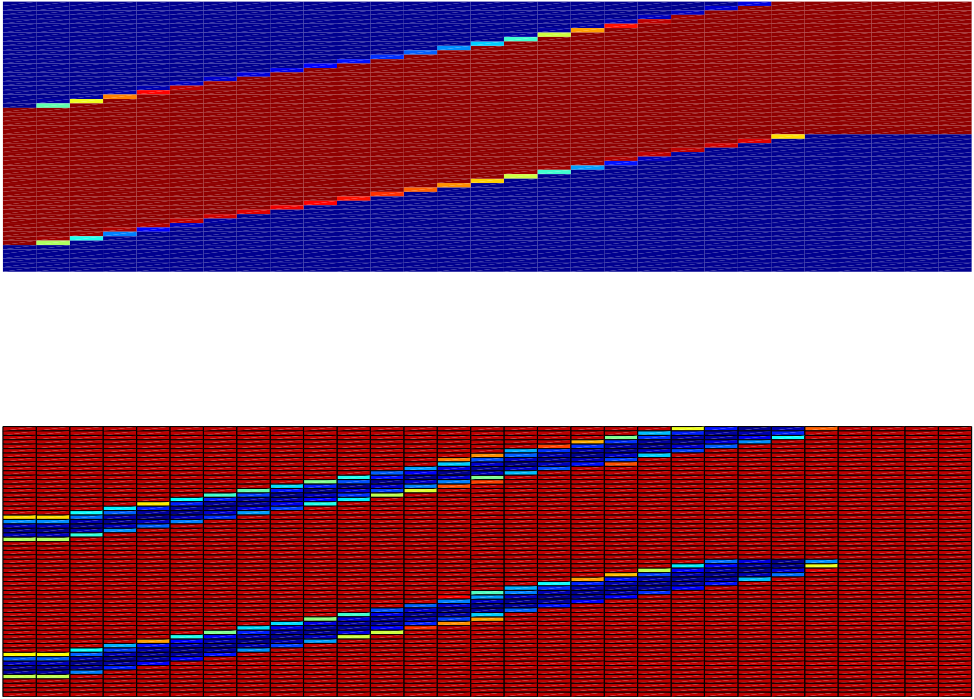}
\caption{A space-time plot of the solution (above) and time-steps
(below) for the propagating front problem, with time going to the
right. The two parts of the plots represent the components for the
two species $A_1$ (lower parts) and $A_2$ (upper parts).}
\label{fig:front-kU}
\end{center}
\end{figure}

\vspace{-1pc}We investigated the passing of a large comet close to
Earth and the Moon and found that the stability factors increase
dramatically at the time $t'$ when the comet comes very close to
the Moon. The conclusion is that if we want to compute accurately
to a point beyond $t'$, we have to be much more careful than if
we want to compute to a point just before $t'$. This is not
evident from the size of residuals or local errors. This is an
example of a Hamiltonian system for which the error growth is
neither linear nor quadratic.

\subsection{A propagating front problem}
The system of PDEs
\begin{equation}
\left\{
\begin{array}{rcl}
\dot{u}_1 - \epsilon u_1'' &=& -u_1 u_2^2, \\
\dot{u}_2 - \epsilon u_2'' &=& u_1 u_2^2
\end{array}
\right.
\end{equation}
on $(0,1)\times(0,T]$ with $\epsilon=0.00001$, $T=100$, and
homogeneous Neumann boundary conditions at $x=0$ and $x=1$
models isothermal autocatalytic reactions (see
\cite{robban:phd}) $A_1 + 2A_2 \rightarrow A_2 + 2A_2$. We
choose the initial conditions as
\begin{displaymath}
u_1(x,0) =
\left\{
\begin{array}{ll}
0, \quad & x<x_0,\\
1, \quad & x\geq x_0,
\end{array}\right.
\end{displaymath}
with $x_0=0.2$ and $u_2(x,0)=1-u_1(x,0)$. An initial reaction
where substance $A_1$ is consumed and substance $A_2$ is formed
will then occur at $x=x_0$, resulting in a decrease in the
concentration $u_1$ and an increase in the concentration $u_2$.
The reaction then propagates to the right until all of substance
$A_1$ is consumed and we have $u_1=0$ and $u_2=1$ in the entire
domain.

\begin{figure}[t]
\begin{center}
\leavevmode
\psfrag{x}{$x$}
\psfrag{U1}{$U_1$}
\psfrag{U2}{$U_2$}
\psfrag{k1}{$k_1$}
\psfrag{k2}{$k_2$}
\includegraphics[width=\textwidth]{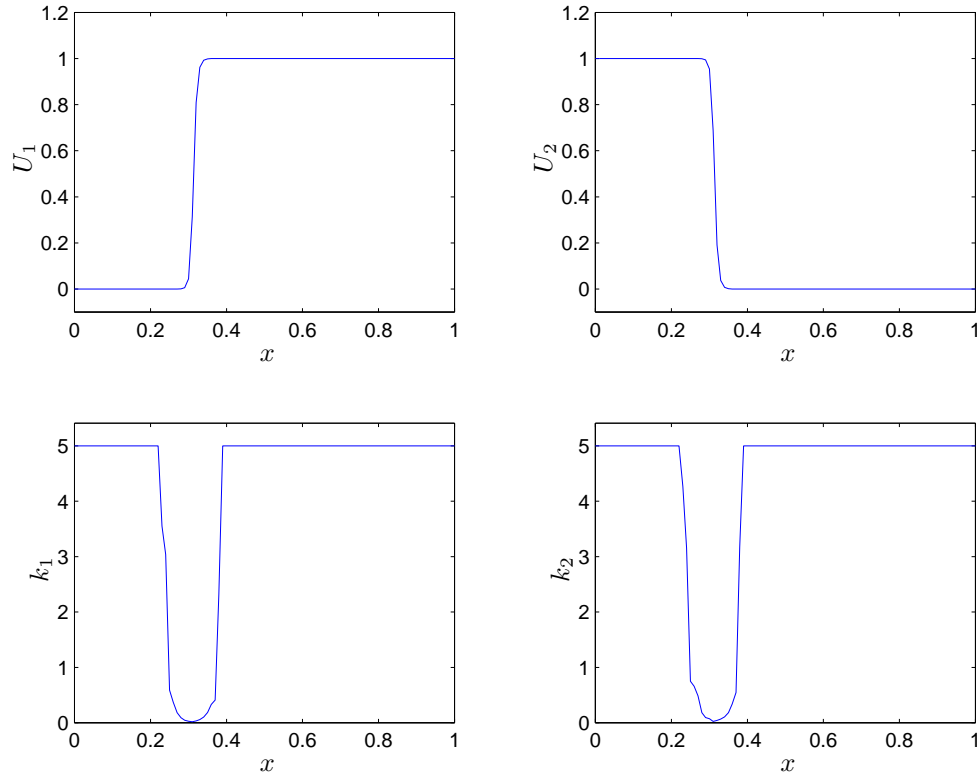}
\caption{The concentrations of the two species
$U_1$ and $U_2$ at time $t=50$ as function of space (above) and
the corresponding time-steps (below).}
\label{fig:front-solution}
\end{center}
\end{figure}

Solving with the $\mathrm{mcG}(2)$ method, we find that the
time-steps are small only close to the reaction front; see
Figures \ref{fig:front-kU} and \ref{fig:front-solution}. The
reaction front propagates to the right as does the domain of
small time-steps.

It is clear that if the reaction front is localized in space and
the domain is large, there is a lot to gain by using small
time-steps only in this area. To verify this, we compute the
solution to within an accuracy of $10^{-7}$ for the final time
error with constant time-steps $k_i(t)=k_0$ for all components
and compare with the multi-adaptive solution. Computing on a
space grid consisting of $16$ nodes on $[0,1]$ (resulting in a
system of ODEs with $32$ components), the solution is computed in
$2.1$ seconds on a regular workstation. Computing on the same
grid with the multi-adaptive method (to within the same
accuracy), we find that the solution is computed in $3.0$
seconds. More work is thus required to compute the multi-adaptive
solution, and the reason is the overhead resulting from
additional bookkeeping and interpolation in the multi-adaptive
computation. However, increasing the size of the domain to $32$
nodes on $[0,2]$ and keeping the same parameters otherwise, we find
the solution is now more localized in the domain and we expect the
multi-adaptive method to perform better compared to a standard method. Indeed, the
computation using equal time-steps now takes $5.7$ seconds,
whereas the multi-adaptive solution is computed in $3.4$ seconds.
In the same way as previously shown in section
\ref{sec:testproblem}, adding extra degrees of freedom does not
substantially increase the cost of solving the problem, since the
main work is done time-stepping the components, which use small
time-steps.

\subsection{Burger's equation with moving nodes}

As a final example, we present a computation in which we combine
multi-adaptivity with the possibility of moving the nodes in a
space discretization of a time-dependent PDE. Solving Burger's
equation,
\begin{equation}
\dot{u} + \mu u u' - \epsilon u'' = 0,
\end{equation}
on $(0,1)\times(0,T]$ with initial condition
\begin{equation}
u_0(x) =
\left\{
\begin{array}{ll}
\sin(\pi x / x_0), & 0 \leq x \leq x_0, \\
0                 & \mbox{elsewhere},
\end{array}
\right.
\end{equation}
and with $\mu=0.1$, $\epsilon=0.001$, and $x_0 = 0.3$, we find the
solution is a shock forming near $x=x_0$. Allowing individual
time-steps within the domain, and moving the nodes of the space
discretization in the direction of the convection, $(1,\mu u)$, we
make the ansatz
\begin{equation}
U(x,t) = \sum_{i=1}^N \xi_i(t) \varphi_i(x,t),
\label{eq:pde,ansatz}
\end{equation}
where the $\{\xi_i\}_{i=1}^N$ are the individual components
computed with the multi-adaptive method, and the
$\{\varphi_i(\cdot,t)\}_{i=1}^N$ are piecewise linear basis
functions in space for any fixed~$t$.

\begin{figure}[t]
\begin{center}
\leavevmode
\psfrag{x}{$x$}
\psfrag{U(x)}{$U(x)$}
\includegraphics[width=\textwidth]{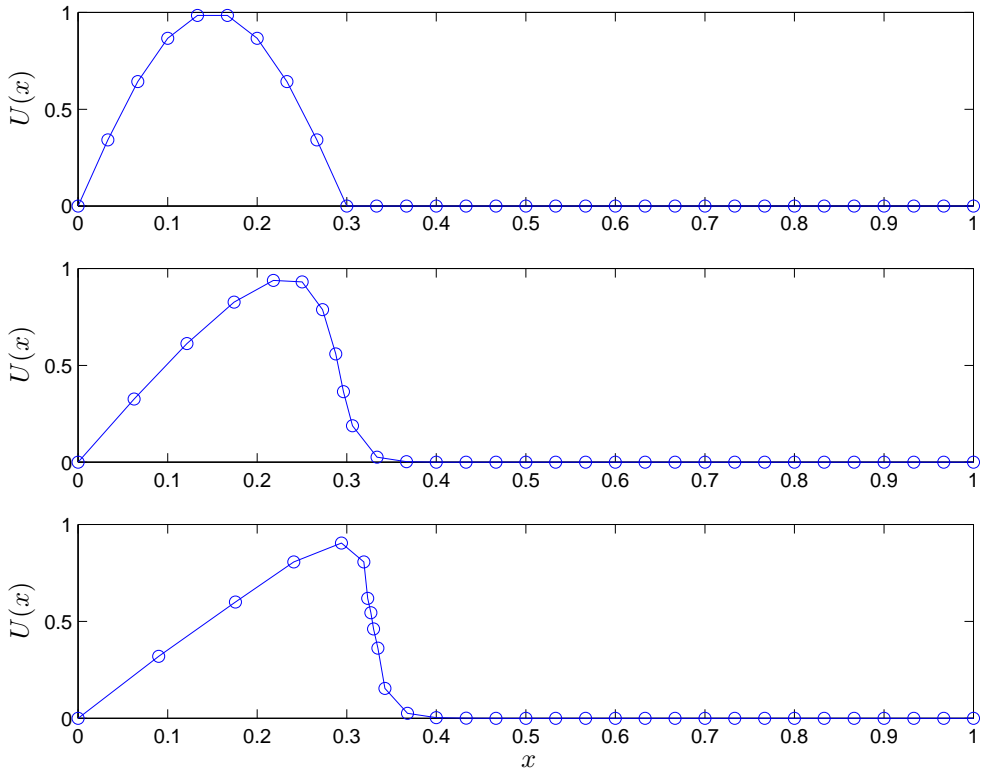}
\caption{The solution to Burger's equation as
function of space at $t=0$, $t=T/2$, and $t=T$.}
\label{fig:burger-solutions}
\end{center}
\end{figure}

\begin{figure}[h!]
\begin{center}
\leavevmode
\psfrag{x}{$x$}
\psfrag{t}{$t$}
\includegraphics[width=\textwidth]{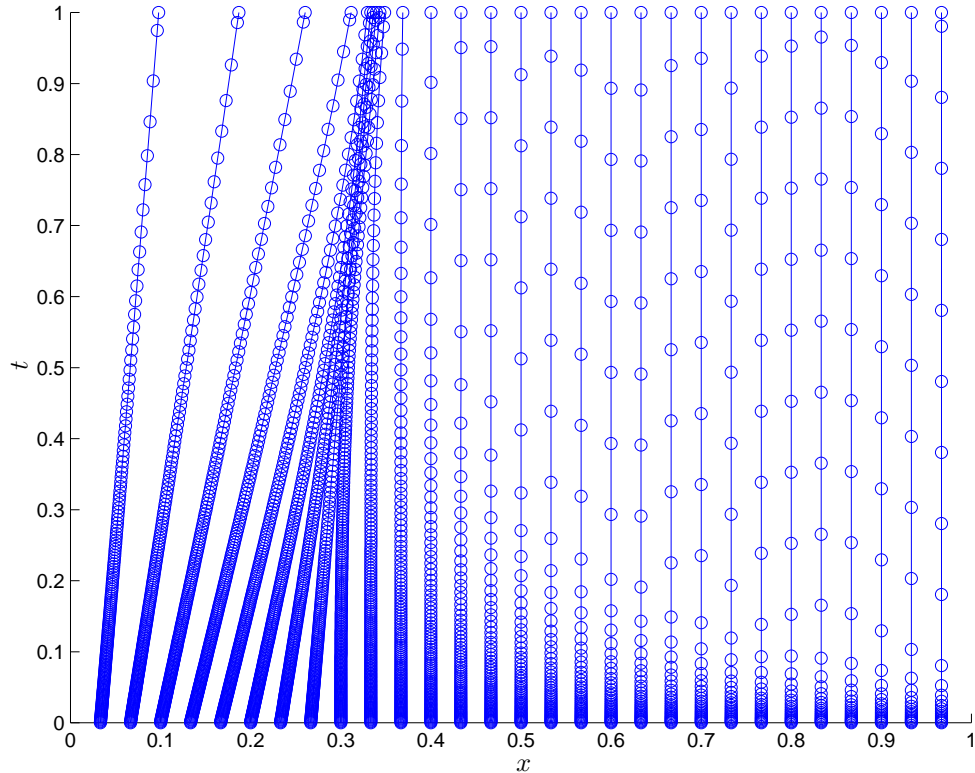}
\caption{Node paths for the multi-adaptive moving-nodes solution
of Burger's equation.}\vspace{-1pc}
\label{fig:burger-paths}
\end{center}
\end{figure}

Solving with the $\mathrm{mdG}(0)$ method, the nodes move into the shock, in the direction
of the convection, so what we are really solving is a heat equation
with multi-adaptive time-steps along the streamlines;
see Figures \ref{fig:burger-solutions} and \ref{fig:burger-paths}.

\section{Future work}
Together with the companion paper \cite{logg:multiadaptivity:I}
(see also
\cite{logg:exjobb,logg:oxford}), this paper serves as a starting
point for further investigation of the multi-adaptive Galerkin
methods and their properties.

Future work will include a more thorough investigation of the
application of the multi-adaptive methods to stiff ODEs, as well
as the construction of efficient multi-adaptive solvers for
time-dependent PDEs, for which memory usage becomes an\linebreak
important issue.


\end{document}